\documentclass[11pt]{amsart}
\usepackage{color}
\usepackage{amsmath}
\usepackage{amssymb}
\usepackage{amsfonts}
\usepackage[latin1]{inputenc}

\DeclareMathSymbol{\subsetneqq}{\mathbin}{AMSb}{36}

\newtheorem{theorem}{Théorème}
\theoremstyle{plain}

\newtheorem{definition}{Définition}

\newtheorem{lemma}{Lemme}
\newtheorem{notation}{Notation}

\newtheorem{proposition}{Proposition}
\newtheorem{remark}{Remarque}

\numberwithin{equation}{section}
\input{tcilatex}

\begin{document}
\title[Extension d'une classe d'unicité]{Extension d'une classe d'unicité
pour les équations de Navier-Stokes}
\author{Ramzi May}
\address{Département de Mathématiques, Faculté des Sciences de Bizerte,
Jarzouna 2001 Bizerte Tunisie.}
\email{Ramzi.May@fsb.rnu.tn}
\date{9 avril 2009}
\maketitle

\noindent \textbf{Résumé:} Récemment, Q. Chen, C. Miao et Z. Zhang \cite{CMZ}
ont montré l'unicité des solutions faibles de Leray dans l'espace $L^{\frac{2%
}{1+r}}([0,T],B_{\infty }^{r,\infty }(\mathbb{R}^{3}))$ avec $r\in ]-\frac{1%
}{2},1].$ Nous proposons dans le présent travail d'étendre ce critère
d'unicité au cas $r\in ]-1,-\frac{1}{2}].\bigskip $

\noindent \textbf{Abstract:} Recently, Q. Chen, C. Miao et Z. Zhang \cite%
{CMZ} have proved that weak Leray solutions of the Navier-Stokes are unique
in the classe $L^{\frac{2}{1+r}}([0,T],B_{\infty }^{r,\infty }(\mathbb{R}%
^{3}))$ with $r\in ]-\frac{1}{2},1].$ In this paper, we establish that this
criterion remains true for $r\in ]-1,-\frac{1}{2}].$

\section{Introduction et énoncé des résultats}

On considère les équations de Navier-Stokes pour un fluide incompressible é%
voluant dans l'espace entier $\mathbb{R}^{d},~d\geq 2,$%
\begin{equation}
\left\{ 
\begin{array}{c}
\partial _{t}u-\Delta u+\nabla (u\otimes u)+\overrightarrow{\nabla }p=0 \\ 
\overrightarrow{\nabla }\cdot u=0 \\ 
u(0,.)=u_{0}(.)%
\end{array}%
\right.  \tag{NS}
\end{equation}%
où $u_{0}$ est la vitesse initiale des particules du fluide, $u=u(t,x)$ dé%
signe la vitesse d'une particule située au point $x\in \mathbb{R}^{d}$ à
l'instant $t\geq 0,$ $p=p(t,x)$ est la pression au point $x\in \mathbb{R}%
^{d} $ à l'instant $t\geq 0,$ $\overrightarrow{\nabla }$ est l'opérateur
gradient, $\overrightarrow{\nabla }\cdot $ est l'opérateur divergence et $%
\nabla (u\otimes u)$ est la fonction vectorielle définie par: $\nabla
(u\otimes u)=(w_{1},\cdots ,w_{d})$ et 
\begin{equation*}
w_{i}\equiv \sum_{k=1}^{d}\frac{\partial }{\partial x_{k}}\left(
u_{k}u_{i}\right) =\overrightarrow{\nabla }\cdot \left( u_{i}u\right) .
\end{equation*}

Rappelons maintenant la notion des solutions faibles des équations de
Navier-Stokes que nous adoptons dans notre travail.

\begin{definition}
Soient $T\in ]0,+\infty ]$ et $u_{0}=(u_{01},\cdots ,u_{0d})$ une
distribution tempérée à divergence nulle. Une solution faible sur $]0,T[$
des équations (NS) est une fonction mesurable $u:Q_{T}\equiv ]0,T[\times $ $%
\mathbb{R}^{d}\rightarrow $ $\mathbb{R}^{d}$ qui vérifie les propriétés
suivantes:

\begin{enumerate}
\item $u\in L_{loc}^{2}(\overline{Q_{T}})$ où $\overline{Q_{T}}\equiv
\lbrack 0,T[\times $ $\mathbb{R}^{d}.$

\item $u\in C([0,T[,S^{\prime }(\mathbb{R}^{d}))$.

\item $u(0)=u_{0}.$

\item $\forall t\in \lbrack 0,T[,~\overrightarrow{\nabla }\cdot u(t)=0$ dans 
$S^{\prime }(\mathbb{R}^{d}).$

\item $\exists p\in D^{\prime }(Q_{T})$ tel que $\partial _{t}u-\Delta
u+\nabla (u\otimes u)+\overrightarrow{\nabla }p=0$ dans $D^{\prime }(Q_{T}).$
\end{enumerate}
\end{definition}

En 1934, J. Leray \cite{Ler} a montré que, pour toute donnée initiale $u_{0}$
appartenant à $L^{2}\left( \mathbb{R}^{d}\right) $ et de divergence nulle,
les équations de Navier-Stokes admettent au moins une solution $u$ faible et
globale en temps qui appartient, pour tout $T>0,$ à l'\textit{espace d'é%
nergie} $\mathcal{L}_{T}$ défini par: 
\begin{equation*}
\mathcal{L}_{T}\equiv L^{\infty }([0,T],L^{2}\left( \mathbb{R}^{d}\right)
)\cap L^{2}([0,T],H^{1}\left( \mathbb{R}^{d}\right) ).
\end{equation*}%
Ceci nous amène à introduire la notion suivante de solutions faibles de
Leary.

\begin{definition}
Soient $T>0$ et $u_{0}\in L^{2}\left( \mathbb{R}^{d}\right) $ telle que $%
\overrightarrow{\nabla }\cdot u_{0}=0.$ On appele une solution faible de
Leray sur $]0,T[$ des équtaions (NS) toute solution faible sur $]0,T[$ des é%
quations (NS) appartenant à l'espace $\mathcal{L}_{T}$.
\end{definition}

Naturellement se pose la question de l'unicité des solutions faibles de
Leray. Dans le cas où la dimension de l'espace $\mathbb{R}^{d}$ est égale à $%
2,$ Il est bien connu que ces solutions sont uniques (voir par exemple \cite%
{Tem}). Dans le cas où la dimension de $\mathbb{R}^{d}$ est supérieure à 3,
la question reste d'une trés grande actualité. On dispose seulement des ré%
ponses partielles. En effet, l'unicité est obtenue sous des hypothèses supplé%
mentaires sur la régularité des solutions. A titre d'exemples, nous citons
les traveaux de J. Serrin \cite{Ser}, de W. Von Wahl \cite{Wah}, de J-Y.
Chemin \cite{Che99} et de I. Gallagher et F. Plonchon \cite{GP}. Dans cette
direction, Q. Chen, C. Miao et Z. Zhang \cite{CMZ} viennent de montrer le ré%
sultat d'unicité suivant.

\begin{theorem}[voir Théorème 1.4 de \protect\cite{CMZ}]
\label{th1} Soient $T>0$ et $u_{0}\in L^{2}\left( \mathbb{R}^{d}\right) $
telle que $\overrightarrow{\nabla }\cdot u_{0}=0.$ Soient $u_{1}$ et $u_{2}$
deux solutions faibles de Leray sur $]0,T[$ des équations (NS). On suppose
que 
\begin{equation*}
u_{1}\in L^{\frac{2}{1-r_{1}}}([0,T],B_{\infty }^{-r_{1},\infty }(\mathbb{R}%
^{3}))\text{ et }u_{2}\in L^{\frac{2}{1-r_{2}}}([0,T],B_{\infty
}^{-r_{2},\infty }(\mathbb{R}^{3}))
\end{equation*}%
où $r_{1},r_{2}\in \lbrack 0,1[$ tels que $r_{1}+r_{2}<1.$ Alors $%
u_{1}=u_{2} $ sur l'intervalle $[0,T].$
\end{theorem}

Il en résulte en particulier que les espaces $L^{\frac{2}{1-r}%
}([0,T],B_{\infty }^{-r,\infty }(\mathbb{R}^{3}))$ avec $r\in \lbrack 0,%
\frac{1}{2}[$ forment une classe d'unicité pour les solutions faibles de
Leray. Dans le présent travail, nous nous proposons d'étendre ce critère
d'unicité aux valeurs $r\in \lbrack \frac{1}{2},1].$ Ce qui répond
positivement à la question posée dans la Remarque 1.7 de \cite{CMZ}.\medskip

Avant d'énoncer nos résultats, nous introduisons la notation suivante:

\begin{notation}
Soient $T>0$ et $r\in ]0,1].$ On pose%
\begin{eqnarray*}
\mathcal{P}_{r,T} &=&L^{\frac{2}{1-r}}([0,T],B_{\infty }^{-r,\infty }(%
\mathbb{R}^{d}))\text{ si }r\neq 1 \\
\mathcal{P}_{r,T} &=&C([0,T],B_{\infty }^{-1,\infty }(\mathbb{R}^{d}))\text{
si }r=1.
\end{eqnarray*}
\end{notation}

\begin{theorem}
\label{th2} On suppose ici que $d\leq 4.$ Soient $T>0$ et $u_{0}\in
L^{2}\left( \mathbb{R}^{d}\right) $ telle que $\overrightarrow{\nabla }\cdot
u_{0}=0.$ Soient $u_{1}$ et $u_{2}$ deux solutions faibles de Leray sur $%
]0,T[$ des équations (NS). On suppose qu'il existe $r_{1},r_{2}\in ]0,1]$
tels que $u_{1}\in \mathcal{P}_{r_{1},T}$ et $u_{2}\in \mathcal{P}%
_{r_{2},T}. $ Alors $u_{1}=u_{2}$ sur l'intervalle $[0,T].$
\end{theorem}

Grâce aux inégalités de Sobolev précisées de Gérard-Meyer et Oru \cite{GMO},
ce théorème va être une conséquence du résultat d'unicité plus général
suivant.

\begin{theorem}
\label{th3} Soient $T>0$, $\left( r_{1},r_{2}\right) \in ]0,1]^{2}$ et $%
(p_{i},q_{i})_{i=1,2}$ deux couples de réels tels que, pour tout $i,$ $%
q_{i}\geq d$ et $p_{i}\geq \frac{4}{1+r_{i}}$ si $r_{i}\neq 1$ et $p_{i}>2$
si $r_{i}=1.$ Soient $u_{1}$ et $u_{2}$ deux solutions faibles sur $]0,T[$
des équations (NS) associées à la même donnée initiale $u_{0}$ telles que%
\begin{equation*}
u_{1}\in L^{p_{1}}([0,T],L^{q_{1}}(\mathbb{R}^{d}))\cap \mathcal{P}_{r_{1},T}%
\text{ et }u_{2}\in L^{p_{2}}([0,T],L^{q_{2}}(\mathbb{R}^{d}))\cap \mathcal{P%
}_{r_{2},T}.
\end{equation*}%
Alors $u_{1}=u_{2}$ sur l'intervalle $[0,T].$
\end{theorem}

La démonstration de ce théorème repose essentiellement sur le résultat de ré%
gularité suivant.

\begin{theorem}
\label{th4} Soient $T>0,r\in ]0,1]$ et $q\geq d$ des réels et soit $p$ un ré%
el supérieur à $\frac{4}{1+r}$ si $r\neq 1$ et strictement supérieur à $2$
si $r=1.$ Si $u\in L^{p}([0,T],L^{q}(\mathbb{R}^{d}))\cap \mathcal{P}_{r,T}$
est une solution faible sur $]0,T[$des équations (NS), alors $\sqrt{t}u\in
L^{\infty }([0,T],L^{\infty }(\mathbb{R}^{d}))$ et $\sqrt{t}\left\Vert
u(t)\right\Vert _{\infty }$ tend vers $0$ lorsque $t\rightarrow 0.$
\end{theorem}

\begin{remark}
Ce théorème implique, en particulier, que toute solution $u$ des équations
de Navier-Stokes appartenant à l'espace $L^{p}([0,T],L^{q}(\mathbb{R}%
^{d}))\cap \mathcal{P}_{r,T}$ est, en effet, une solution classique des é%
quations (NS) sur le cylindre $\overline{Q_{T}}\equiv ]0,T]\times \mathbb{R}%
^{d}$ i.e $u\in C^{\infty }(\overline{Q_{T}})$ (voir la démonstration de ce
théorème dans la dernière section de cet article).
\end{remark}

\begin{remark}
\label{Rq}Dans le cas où $r=1,$ le théorème \ref{th4} a été récemment démontr%
é par P.G. Lemarié-Rieusset \cite{Lem07} lorsque $q>d$ puis par l'auteur de
cet article \cite{May09} dans le cas limite $q=d.$ Ainsi, on se limitera,
dans la suite de ce travail, à démontrer ce théorème dans le cas $r\in
]0,1[. $ Nous signalons que la démonstration du théoréme dans ce cas utilse
des techniques différentes de celles utilisées dans \cite{Lem07} et dans 
\cite{May09} pour montrer le théorème dans le cas limite $r=1.$
\end{remark}

Cet article s'organise comme suit: dans la section suivante, nous rappelons
en premier lieu la notion des solutions mild des équations de Navier-Stokes
introduite dans \cite{FLT}, puis nous introduisons les espaces de Besov et
les espaces de Chemin-Lerner et nous citons quelques propriétés de ces
espaces dont on fait usage. Dans la troisième section, nous montrons comment
le théorème principal \ref{th4} implique le théorème \ref{th3} puis le théorè%
me \ref{th2}. La dernière section est consacrée à la preuve du théorème \ref%
{th4} dans le cas où $r\in ]0,1[.$

\section{Préliminaires}

\subsection{Notations}

\begin{enumerate}
\item Tous les espaces fonctionels utilisés dans ce travail sont définis sur
l'espace entier $\mathbb{R}^{d}.$ Ainsi, pour alléger l'écriture, nous é%
crivons par exemple $L_{x}^{q},H^{s}$ et $B_{p}^{s,q}$ pour désigner
respectivement les espaces $L^{q}\left( \mathbb{R}^{d}\right) ,H^{s}\left( 
\mathbb{R}^{d}\right) $ et $B_{p}^{s,q}\left( \mathbb{R}^{d}\right) .$

\item Si $X$ est un espace vectoriel et $n\in \mathbb{N}^{\ast },$ on écrit
souvent $(f_{1},\cdots ,f_{n})\in X$ à la place de $(f_{1},\cdots ,f_{n})\in
X^{n}.$

\item Si $X$ est un espace de Banach, $T$ est un réel positif et $p\in
\lbrack 1,+\infty ],$ on notera par $L_{T}^{p}(X)$ et $L_{T}^{p}X$ l'espace $%
L^{p}([0,T],X).$

\item Soit $p$ un réel supérieur à $1$. On désigne par $\mathbf{E}_{p}$
l'espace de fonctions $f\in L_{loc}^{p}(\mathbb{R}^{d})$ telles que:%
\begin{equation*}
\left\Vert f\right\Vert _{\mathbf{E}_{p}}\equiv \sup_{x_{0}\in \mathbb{R}%
^{d}}\left\Vert 1_{B(x_{0},1)}f\right\Vert _{p}<\infty \text{ et }%
\lim_{\left\Vert x_{0}\right\Vert \rightarrow \infty }\left\Vert
1_{B(x_{0},1)}f\right\Vert _{p}=0.
\end{equation*}

\item Si $A$ et $B$ sont des fonctions réelles, on écrit $A\lesssim B$
lorsque $A\leq CB$ où $C$ une constante strictement positive indépendante
des paramètres qui définissent $A$ et $B$.
\end{enumerate}

\subsection{Solutions mild des équations de Navier-Stokes}

On désigne par $\mathbb{P}$ le projecteur de Leray sur l'espace des
distributions à divergence nulle. On rappelle que $\mathbb{P=}\left( \mathbb{%
P}_{ij}\right) _{1\leq i,j\leq d}$ est défini à l'aide des transformations
de Riesz $\left( \mathcal{R}_{i}\right) _{1\leq i\leq d}$ par la relation:%
\begin{equation*}
\mathbb{P}_{ij}\left( f\right) =\delta _{ij}f-\mathcal{R}_{i}\mathcal{R}%
_{j}\left( f\right)
\end{equation*}%
où $\delta _{ij}$ est le symbole de Kronecker.

Soit $u_{0}=(u_{01},\cdots ,u_{0d})\in S^{\prime }(\mathbb{R}^{d})$ une
distribution à divergence nulle. En appliquant formellement l'opérateur de
Leray $\mathbb{P}$ aux équations (NS) on obtient le système d'équations
suivant:

\begin{equation*}
\left\{ 
\begin{array}{c}
\partial _{t}u-\Delta u=-\mathbb{P}\nabla (u\otimes u) \\ 
u(0,.)=u_{0}(.).%
\end{array}%
\right.
\end{equation*}%
Ensuite, en utilisant la formule de Duhamel, nous tronsformons formellement
ces équations en des équations intégrales:%
\begin{equation}
u(t)=e^{t\Delta }u_{0}+\mathbf{B}(u,u)(t)  \tag{NSI}
\end{equation}%
où $\left( e^{t\Delta }\right) _{t\geq 0}$ est le semi-groupe de la chaleur
et $\mathbf{B}$ est l'application bilinéaire définie par: 
\begin{equation*}
\mathbf{B}(u,v)(t)=\mathbb{L}_{Oss}(u\otimes v)(t)
\end{equation*}%
$\mathbb{L}_{Oss}$ est l'opérateur linéaire défini par:%
\begin{equation}
\mathbb{L}_{Oss}(f)=-\int_{0}^{t}e^{(t-s)\Delta }\mathbb{P}\nabla (f)ds
\label{oss}
\end{equation}%
Dans la suite, on appellera $\mathbb{L}_{Oss}$ \textquotedblleft \textit{l'op%
érateur intégral d'Oseen\textquotedblright }.\medskip

Dans \cite{FLT}, J. Furioli, P.-G. Lemarié-Rieusset et E. Terraneo ont montré
que, dans la classe des solutions $L_{loc}^{2}([0,T[,\mathbf{E}_{2}),$ les é%
quations (NS) et leur forme intégrale (NSI) sont équivalentes$.$ Ce qui nous
conduit à introduire la notion suivante des solutions \textit{mild} que nous
adoptons dans notre travail.

\begin{definition}
Soient $T>0$ et $u_{0}\in S^{\prime }\left( \mathbb{R}^{d}\right) .$ Une
solution mild des équations de Navier-Stokes (NS) sur l'intervalle $]0,T[ $
est une fonction $u\in L_{los}^{2}([0,T[,\mathbf{E}_{2})$ solution dans $%
D^{\prime }(Q_{T})$ des équations intégrales (NSI) où $Q_{T}=]0,T[\times 
\mathbb{R}^{d}.$
\end{definition}

\begin{remark}
\label{rq3}Il est bien connu (voir par exemple \cite{FLT} et \cite{Lem02})
que toute solution mild $u$ des équations de Navier-Stokes sur $]0,T[$
appartient à l'espace $C([0,T[,B_{\infty }^{-d-1,\infty }\left( \mathbb{R}%
^{d}\right) )$ et que, pour tout $t\in \lbrack 0,T[,$%
\begin{equation*}
u(t)=e^{t\Delta }u_{0}+\mathbf{B}(u,u)(t)\text{ dans }B_{\infty
}^{-d-1,\infty }\left( \mathbb{R}^{d}\right) .
\end{equation*}
\end{remark}

\begin{remark}
Toutes les solutions des équations de Navier-Stokes que nous considèrons
dans ce travail sont des solutions mild. Ainsi, dans la suite, étant donnée $%
u\in L_{los}^{2}([0,T[,\mathbf{E}_{2})$, nous écrivons souvent $u$ est une
solution des équations (NS) pour signifier que $u$ est une solution mild sur 
$]0,T[$des équations (NS).
\end{remark}

\begin{remark}
\label{mild}Soit $u$ une solution mild sur $]0,T[$ des équations (NS)$.$ En
utilisant la propriété du semi-groupe de $\left( e^{t\Delta }\right) _{t\geq
0}$, on vérifie aisément que pour tous $t_{0}\leq t\in \lbrack 0,T[,$%
\begin{equation*}
u(t)=e^{(t-t_{0})}u(t_{0})-\int_{t_{0}}^{t}e^{(t-s)\Delta }\mathbb{P}\nabla
(u\otimes u)ds.
\end{equation*}%
Ce qui implique que la fonction $u_{t_{0}}\equiv u(.+t_{0})$ est une
solution mild sur $]0,T-t_{0}[$ des équations de Navier-Stokes associées à
la donnée initiale $u(t_{0}).$
\end{remark}

\begin{remark}
Dans la suite de notre travail, nous oublierons la condition de nullité du
champs de vecteur $u$ qui n'aura aucun rôle.
\end{remark}

\subsection{Espaces de Besov et Espaces de Chemin-Lerner}

Nous rappelons tout d'abord la décomposition de Littlewood-Paley. Soit $%
\varphi \in C_{c}^{\infty }(\mathbb{R}^{d})$ égale à 1 sur un voisinage de
l'origine$.$ On considère ensuite la fonction $\psi \in C_{c}^{\infty }(%
\mathbb{R}^{d})$ définie par $\psi (\xi )=\varphi (\frac{\xi }{2})-\varphi
(\xi ).$ Pour $j\in \mathbb{N},$ on désigne par $S_{j}$ et $\Delta _{j}$ les
multiplicateurs de Fourier définis pour $f\in S^{\prime }(\mathbb{R}^{d})$
et $v\in S^{\prime }(\mathbb{R}\times \mathbb{R}^{d})$ par:%
\begin{eqnarray*}
S_{j}f &=&\mathcal{F}^{-1}\left( \varphi (\frac{\xi }{2^{j}})\mathcal{F}%
(f)\right) ,~\Delta _{j}f=\mathcal{F}^{-1}\left( \psi (\frac{\xi }{2^{j}})%
\mathcal{F}(f)\right) \\
S_{j}v &=&\mathcal{F}_{x}^{-1}\left( \varphi (\frac{\xi }{2^{j}})\mathcal{F}%
_{x}(v)\right) ,~\Delta _{j}v=\mathcal{F}_{x}^{-1}\left( \psi (\frac{\xi }{%
2^{j}})\mathcal{F}_{x}(v)\right)
\end{eqnarray*}%
où $\mathcal{F}$ et $\mathcal{F}^{-1}$ sont la transformation de Fourier et
son inverse définies sur $S^{\prime }(\mathbb{R}^{d})$ et $\mathcal{F}_{x}$
et $\mathcal{F}_{x}^{-1}$ sont la transformation de Fourier et son inverse
par rapport à la variable $x$ définies sur $S^{\prime }(\mathbb{R}\times 
\mathbb{R}^{d})$.

\begin{notation}
Dans la suite, on notera souvent l'opérateur $S_{0}$ par $\Delta _{-1}.$
\end{notation}

Nous pouvons maintenant rappeler la définition d'une classe des espaces de
Besov.

\begin{definition}
Soient $s\in \mathbb{R}$ et $q\in \lbrack 1,+\infty ].$ L'espace de Besov $%
B_{q}^{s,\infty }$ est l'espace de distributions $f\in S^{\prime }(\mathbb{R}%
^{d})$ telles que%
\begin{equation*}
\left\Vert f\right\Vert _{B_{q}^{s,\infty }}\equiv \sup_{j\geq
-1}2^{sj}\left\Vert \Delta _{j}f\right\Vert _{q}<\infty .
\end{equation*}%
On désigne par $\tilde{B}_{q}^{s,\infty }$ l'adherence de $S(\mathbb{R}^{d})$
dans $B_{q}^{s,\infty }.$
\end{definition}

Nous rappelons aussi la définition d'une classe des espaces de Chemin-Lerner
(\cite{Dan}, \cite{Che99}, \cite{CL}).

\begin{definition}
Soient $T>0,s\in \mathbb{R}$ et $p,q\in \lbrack 1,+\infty ].$ L'espace de
Chemin-Lerner $\tilde{L}_{T}^{p}B_{q}^{s,\infty }$ est l'espace de
distributions $v\in S^{\prime }(\mathbb{R}\times \mathbb{R}^{d})$ telles que%
\begin{equation*}
\left\Vert v\right\Vert _{\tilde{L}_{T}^{p}B_{q}^{s,\infty }}\equiv
\sup_{j\geq -1}2^{sj}\left\Vert \Delta _{j}v\right\Vert
_{L_{T}^{p}L_{x}^{q}}<\infty .
\end{equation*}%
On désigne par $\mathbf{\tilde{L}}_{T}^{p}B_{q}^{s,\infty }$ l'ensemble de
distributions $v\in \tilde{L}_{T}^{p}B_{q}^{s,\infty }$ telles que:%
\begin{equation*}
\left\Vert v\right\Vert _{\tilde{L}_{T_{1}}^{p}B_{q}^{s,\infty }}\rightarrow
0\text{ lorsque }T_{1}\rightarrow 0.
\end{equation*}
\end{definition}

La proposition suivante regroupe quelques propriétés simples et utiles des
espace de Besov et de Chemin-Lerner.

\begin{proposition}
\label{pro1} Soient $T>0,s\in \mathbb{R},(p,q)\in \lbrack 1,+\infty ]^{2}$
et $p_{1}\in \lbrack 1,+\infty \lbrack .$ Les assertions suivantes sont
vraies:

\begin{enumerate}
\item $L_{T}^{p}B_{q}^{s,\infty }\subset \tilde{L}_{T}^{p}B_{q}^{s,\infty
},L_{T}^{\infty }B_{q}^{s,\infty }=\tilde{L}_{T}^{\infty }B_{q}^{s,\infty }$
et $L_{T}^{p_{1}}B_{q}^{s,\infty }\subset \mathbf{\tilde{L}}%
_{T}^{p_{1}}B_{q}^{s,\infty }.$

\item Les opérateurs $\mathbb{P}_{ij}\frac{\partial }{\partial x_{k}}$
envoient continûment $B_{q}^{s,\infty }$ (respectivement $%
L_{T}^{p}B_{q}^{s,\infty }$) dans $B_{q}^{s-1,\infty }$ (respectivement $%
L_{T}^{p}B_{q}^{s-1,\infty }$).

\item \lbrack Injections de Bernstein] Pour tout $m\in \lbrack q,\infty ],$
on a 
\begin{equation*}
B_{q}^{s,\infty }\subset B_{m}^{s+d(\frac{1}{m}-\frac{1}{q}),\infty }\text{
et }\tilde{L}_{T}^{p}B_{q}^{s,\infty }\subset \tilde{L}_{T}^{p}B_{m}^{s+d(%
\frac{1}{m}-\frac{1}{q}),\infty }
\end{equation*}%
En plus, les normes des injections et des applications linéaires considérées
dans cette proposition ne dépendent pas de $T$.
\end{enumerate}
\end{proposition}

La démonstration de cette proposition est tout à fait classique.\medskip

Il est bien connu (voir par exemple \cite{Can}, \cite{Lem02} et \cite{Tri})
que le semi-groupe de la chaleur $\left( e^{t\Delta }\right) _{t\geq 0}$
permet de caractériser les espaces de Besov. Nous rappelons dans la
proposition suivante un cas particulier de cette caractérisation.

\begin{proposition}
\label{pro2} Soient $q\in \lbrack 1,+\infty ]$ et $s$ un réel strictement
positif. Alors, pour tout $\delta >0,$ la quantité 
\begin{equation*}
\sup_{0<\theta <\delta }\theta ^{\frac{s}{2}}\left\Vert e^{\theta \Delta
}f\right\Vert _{q}
\end{equation*}%
définit une norme sur $B_{q}^{-s,\infty }$ équivalente à la norme originale.
\end{proposition}

Pour pouvoir établir des estimations de produit de type "estimations douces"
dans les espaces de Besov et les espaces de Chemin-Lerner, nous introduisons
une version simplifiée du paraproduit de Bony: Pour $f$ et $g$ dans $%
S^{\prime }(\mathbb{R}^{d})$ (ou dans $S^{\prime }(\mathbb{R}\times \mathbb{R%
}^{d})$)$,$ on définit formellement $\Pi _{1}(f,g)$ et $\Pi _{2}(f,g)$ par:%
\begin{equation*}
\Pi _{1}(f,g)=\sum_{j\geq -1}S_{j+1}f\Delta _{j}g\text{ et }\Pi
_{2}(f,g)=\sum_{j\geq 0}S_{j}f\Delta _{j}g.
\end{equation*}%
On obtient, au moins formellement, l'égalité suivante: $fg=\Pi _{1}(f,g)+\Pi
_{2}(g,f).$ On appellera dans la suite $\Pi _{1}$ et $\Pi _{2}$ les
\textquotedblleft \textit{opérateurs du paraproduit de
Bony\textquotedblright }.\medskip

La proposition suivante décrit quelques propriétés d'opérances de ces opé%
rateurs sur les espaces de Besov et les espaces de Chemin-Lerner.

\begin{proposition}[Lois du praproduit de Bony]
\label{pro3} Soient $T>0,$ $\sigma _{1\text{ }}<\sigma _{2}$ deux réels
strictement positifs et $(p_{1},q_{1}),(p_{2},q_{2})$ $\in \lbrack 1,+\infty
]^{2}$ tels que $\frac{1}{p}\equiv \frac{1}{p_{1}}+\frac{1}{p_{2}}$ et $%
\frac{1}{q}\equiv \frac{1}{q_{1}}+\frac{1}{q_{2}}$ soient inférieurs à 1.\
Alors on a les deux assertions suivantes:

\begin{enumerate}
\item Les opérateurs du paraproduit de Bony $\Pi _{1}$ et $\Pi _{2}$ sont
continus de $B_{q_{1}}^{-\sigma _{1},\infty }\times B_{q_{2}}^{\sigma
_{2},\infty }$dans l'espace $B_{q}^{\sigma _{2}-\sigma _{1},\infty }.$

\item Les opérateurs du paraproduit de Bony $\Pi _{1}$ et $\Pi _{2}$ sont
continus de $\tilde{L}_{T}^{p_{1}}B_{q_{1}}^{-\sigma _{1},\infty }\times 
\tilde{L}_{T}^{p_{2}}B_{q_{2}}^{\sigma _{2},\infty }$ dans l'espace $\tilde{L%
}_{T}^{p}B_{q}^{\sigma _{2}-\sigma _{1},\infty }$ et de $%
L_{T}^{p_{1}}L_{x}^{q_{1}}\times \tilde{L}_{T}^{p_{2}}B_{q_{2}}^{\sigma
_{2},\infty }$ dans l'espace $\tilde{L}_{T}^{p}B_{q}^{\sigma _{2},\infty }$.
En plus, leurs normes sont indépendantes de $T.$
\end{enumerate}
\end{proposition}

La démonstration de cette proposition est simple, voir par exemple \cite%
{Che99} où des résulats similaires sont prouvés.\medskip

Nous étudions maintenant l'effet régularisant de l'équation de la chaleur
sur les espaces de Besov et les espaces de Chemin-Lerner.\medskip

La première proposition concerne l'effet régularisant du semi-groupe de la
chaleur $(e^{t\Delta })_{t\geq 0}.$

\begin{proposition}[Effet régularisant du semi-groupe de la chaleur]
\label{pro4} Soient $T>0,\left( s,s_{1},s_{2}\right) \in \mathbb{R}^{3}$ et $%
(p,q)\in \lbrack 1,+\infty ]^{2}.$ On a les assertions suivantes:

\begin{enumerate}
\item Si $s_{1}\leq s_{2}$ alors la famille $(t^{\frac{s_{2}-s_{1}}{2}%
}e^{t\Delta })_{0<t\leq T}$ est borneé dans $\mathcal{L}(B_{q}^{s_{1},\infty
},B_{q}^{s_{2},\infty }).$

\item L'opéreteur linéaire $e^{t\Delta }$ envoie continûment $%
B_{q}^{s,\infty }$ dans $\tilde{L}_{T}^{p}B_{q}^{s+\frac{2}{p},\infty }.$ En
plus, si $p<\infty ,$ alors $e^{t\Delta }$ envoie continûment $\tilde{B}%
_{q}^{s,\infty }$ dans $\mathbf{\tilde{L}}_{T}^{p}B_{q}^{s+\frac{2}{p}%
,\infty }.$
\end{enumerate}
\end{proposition}

La seconde proposition concerne l'effet régularisant de l'opérateur intégral
d'Oseen $\mathbb{L}_{oss}$ défini par (\ref{oss}).

\begin{proposition}[Effet régularisant de l'opérateur $\mathbb{L}_{Oss}$]
\label{pro5} Soient $T>0,s\in \mathbb{R}$ et $(p_{1},p_{2},q)\in \lbrack
1,+\infty ]^{3}$ tel que $p_{1}\leq p_{2}.$ On pose $s^{\prime }=s+1-2(\frac{%
1}{p_{1}}-\frac{1}{p_{2}}).$ Alors l'opérateur intégral d'Oseen $\mathbb{L}%
_{oss}$ envoie continûment $\tilde{L}_{T}^{p_{1}}B_{q}^{s,\infty }$ dans $%
\tilde{L}_{T}^{p_{2}}B_{q}^{s^{\prime },\infty }$ et sa norme est majorée
par $C(1+T)$ où $C$ est une constante positive indépendante de $T.$
\end{proposition}

Pour la démonstration de ces deux propositions, nous renvoyons le lecteur
aux références \cite{Che99} et \cite{Dan}.

\section{Démonstration des Théorèmes \protect\ref{th2} et \protect\ref{th3}.}

Dans cette courte section, nous allons voir comment le Théorème \ref{th4}
(qui sera démontré dans la section suivante) permet de montrer le Théorème %
\ref{th3}, puis comment ce dernier implique à son tour le Théorème \ref{th2}.

\subsection{Démonstration du théorème \protect\ref{th3}}

On pose $p=\inf (p_{1},p_{2}).$ D'après le Théorème \ref{th4}, on a pour $%
i=1 $ ou $2,$ 
\begin{equation*}
u_{i}\in L_{T}^{p}\left( \mathbf{E}_{d}\right) ,~\sqrt{t}u_{i}\in
L_{T}^{\infty }L_{x}^{\infty }\text{ et }\lim_{t\rightarrow 0}\sqrt{t}%
\left\Vert u_{i}(t)\right\Vert _{\infty }=0
\end{equation*}%
(voir la sous section 2.1, pour la définition de $\mathbf{E}_{d}$). Ensuite,
en remarquant que $p>2$ et en utilisant le fait que la multiplication par
une fonction $L^{\infty }(\mathbb{R}^{d})$ et la convolution avec une
fonction $L^{1}(\mathbb{R}^{d})$ sont continues sur l'espace $\mathbf{E}%
_{d}, $ on peut facilement adapter la démonstration du Lemme 10 de \cite%
{Lem07} (voir aussi la démonstration du Lemme 5.2 dans \cite{May09}) à notre
cas pour conclure que $u_{1}=u_{2}$ sur $[0,T].$ $\blacksquare $

\subsection{Démonstration du théorème \protect\ref{th2}}

Soit $i=1$ ou $2.$ Tout d'abord, les inégalités classiques d'interpolation
dans les espaces de Lebesgue et dans les espaces de Sobolev impliquent que $%
u_{i}$ appartient à l'espace $L_{T}^{\frac{2}{r_{i}}}H^{r_{i}}.$ Ainsi, en
appliquant la version nonhomogènne des inégalités de Sobolev précisées dues à
P. Gérard, Y. Meyer et F. Oru \cite{GMO} (voir aussi \cite{Lem07}, pour une
autre démonstration des ces inégalités): 
\begin{eqnarray*}
\left\Vert f\right\Vert _{q} &\lesssim &\left( \left\Vert f\right\Vert
_{W^{\alpha ,p}}\right) ^{1-\frac{\alpha }{\beta }}\left( \left\Vert
f\right\Vert _{B_{\infty }^{\alpha -\beta ,\infty }}\right) ^{\frac{\alpha }{%
\beta }} \\
0 &<&\alpha <\beta ,~1<p<\infty ,~\frac{1}{q}=(1-\frac{\alpha }{\beta })p,
\end{eqnarray*}%
avec $\alpha =r_{i},\beta =2r_{i}$ et $p=2,$ et en utilisant ensuite les iné%
galités de Hölder, on trouve que: 
\begin{equation*}
u_{i}\in L^{4}([0,T],L^{4}(\mathbb{R}^{d})).
\end{equation*}%
Ce qui termine la preuve grâce au théorème \ref{th3}. $\blacksquare $

\section{Démonstration du Théorème \protect\ref{th4}}

Nous consacrons cette section à la preuve du Théorème principal \ref{th4}
dans le cas où $r\in ]0,1[$ (voir la Remarque \ref{Rq}).

La démonstration de ce théorème nécessite que l'on établisse auparavant
quelques résultats intermédiaires par ailleurs utiles en eux mêmes.\medskip

La première proposition est un résultat d'unicité locale sous une hypothése
supplémentaire sur la régularité de la donnée initiale.

\begin{proposition}
\label{prop1} Soient $T>0,r\in ]0,1[,q\geq d$ et $p\geq \frac{4}{1+r}$ des ré%
els. Si la donnée initiale $u_{0}\in L^{q}(\mathbb{R}^{d})$ et si $%
u_{1},u_{2}\in L_{T}^{\frac{2}{1-r}}\left( B_{\infty }^{-r,\infty }\right)
\cap L_{T}^{p}L_{x}^{q}$ sont deux solutions mild sur $]0,T[$ des équations
(NS), alors il existe $\delta \in ]0,T]$ tel que $u_{1}=u_{2}$ sur $%
[0,\delta ].$
\end{proposition}

Dans la seconde proposition, nous montrons un résultat de persistance de la r%
égularité de la donnée initiale et un critère du contrôle d'explosion en
temps fini des solutions régulières des équations de Navier-Stokes.

\begin{proposition}
\label{prop2} Soient $q\geq d$ un réel et $u_{0}\in L^{q}(\mathbb{R}^{d}).$
Alors on a les assertions suivantes:

\begin{enumerate}
\item Les équations de Navier-Stokes admettent une et une seule solution
maximale $u\in C([0,T^{\ast }[,L^{q}(\mathbb{R}^{d})).$ En plus, pour tout $%
\sigma >0,$ 
\begin{equation}
u\in C^{\infty }(]0,T^{\ast }[,\tilde{B}_{\infty }^{\sigma ,\infty }).
\label{reg}
\end{equation}

\item Si on suppose en plus que $u_{0}\in B_{\infty }^{-r,\infty }$ avec $%
r\in ]0,1[,$ alors la solution maximale $u$ appartient aussi à l'espace $%
L_{loc}^{\infty }([0,T^{\ast }[,B_{\infty }^{-r.,\infty }).$

\item Si le temps maximal d'existence $T^{\ast }$ est fini$,$ alors, pour
tout $r\in ]0,1[,$ il existe une constante $\varepsilon _{r,d}>0$, qui ne dé%
pend que de $r$ et $d,$ telle que 
\begin{equation}
\underline{\lim }_{t\rightarrow T^{\ast }}(T^{\ast }-t)^{\frac{1-r}{2}%
}\left\Vert u(t)\right\Vert _{B_{\infty }^{-r,\infty }}\geq \varepsilon
_{r,d}.  \label{giga}
\end{equation}%
En particulier, 
\begin{equation*}
\int_{T^{\ast }/2}^{T^{\ast }}\left( \left\Vert u(t)\right\Vert _{B_{\infty
}^{-r,\infty }}\right) ^{\frac{2}{1-r}}dt=+\infty .
\end{equation*}
\end{enumerate}

\begin{remark}
L'estimation (\ref{giga}) améliore un résultat similaire dû à Y. Giga \cite%
{Gig} où la norme dans l'espace de Besov $B_{\infty }^{-r,\infty }$ est
remplacée par la norme dans l'espace $L^{\frac{d}{r}}(\mathbb{R}^{d}).$
\end{remark}
\end{proposition}

La dernière proposition concerne le comportement au voisinage de l'instant
initiale des solutions régulières des équations de Navier-Stokes appartenant 
à l'espace $L_{T}^{\frac{2}{1-r}}B_{\infty }^{-r,\infty }.$

\begin{proposition}
\label{prop3} Soit $r\in ]0,1[,T>0$ et $u\in C(]0,T],B_{\infty }^{1,\infty
})\cap L^{\frac{2}{1-r}}([0,T],B_{\infty }^{-r,\infty })$ une solution mild
sur $]0,T[$ des équations de Navier-Stokes. Alors $\sqrt{t}\left\Vert
u(t)\right\Vert _{\infty }\rightarrow 0$ lorsque $t\rightarrow 0.$
\end{proposition}

Admettons pour le moment ces propositions et voyons comment elles impliquent
le théorème principal.\bigskip

\noindent \textbf{Démonstration du Théorème \ref{th4}:} On pose $\Omega
_{q,r}\equiv \{t_{0}\in ]0,T[$ tel que $u(t_{0})\in L^{q}(\mathbb{R}%
^{d})\cap B_{\infty }^{-r,\infty }\}.$ Soit $t_{0}$ un élement quelconque
mais fixe de $\Omega _{r,q}$. La Proposition \ref{prop2} nous apprend que
les équations de Navier-Stokes associées à la donnée initiale $v_{0}\equiv
u(t_{0})$ admettent une solution maximale $v$ qui appartient à l'espace $%
C([0,T^{\ast }[,L^{q})\cap $ $L_{loc}^{\infty }([0,T^{\ast }[,B_{\infty
}^{-r,\infty })$ où $T^{\ast }$ est son temps maximal d'existence. La
Remarque \ref{mild} et la Proposition \ref{prop1} impliquent alors qu'il
existe $\delta \in ]0,\delta _{0}[$ tel que $v=u(.+t_{0})$ sur $[0,\delta ]$
où $\delta _{0}\equiv \inf \{T^{\ast },T-t_{0}\}.$ Ce qui justife la dé%
finition:%
\begin{equation*}
\delta _{\ast }\equiv \sup \{\delta \in ]0,\delta _{0}[\mid v=u(.+t_{0})%
\text{ sur }[0,\delta ]\}.
\end{equation*}%
Supposons par l'absurde que $\delta _{\ast }<\delta _{0}.$ En utilisant la
continuité de $v$ sur $[0,\delta _{0}[$ à valeurs dans $L^{q}(\mathbb{R}%
^{d}) $ et celle de $u(.+t_{0})$ sur $[0,\delta _{0}[$ à valeurs dans
l'espace de Besov $B_{\infty }^{-1-d,\infty }$ (voir la Remarque \ref{rq3}),
on en déduit alors que $v(\delta _{\ast })=u(\delta _{\ast }+t_{0})\in L^{q}(%
\mathbb{R}^{d}).$ Ainsi, en appliquant à nouveau la Proposition \ref{prop1}
aux équations de Navier Stokes associées à la nouvelle donnée initiale $%
v(\delta _{\ast })$, on trouve qu'il existe $\delta ^{\prime }>\delta _{\ast
}$ tel que $v=u(.+t_{0})$ sur $[0,\delta ^{\prime }].$ Ce qui contredit la dé%
finition de $\delta _{\ast }.$ On conclut alors que $v=u(.+t_{0})$ sur $%
[0,\delta _{0}[.$ Ce qui implique, vu l'hypothèse sur la solution $u,$ que $%
v\in L^{\frac{2}{1-r}}([0,\delta _{0}[,B_{\infty }^{-r,\infty }).$ Par consé%
quent, la troisième assertion de la Proposition \ref{prop2}, implique que $%
T^{\ast }>\delta _{0}$. Il en résulte ainsi que $u(.+t_{0})=v$ sur $%
[0,T-t_{0}].$ Utilisons maintenant la régularité de la solution $v,$ assurée
par la première assertion de la Proposition \ref{prop2}, et la densité de $%
\Omega _{r,q}$ dans $]0,T],$ on conclut que la solution $u$ appartient à $%
\cap _{\sigma >0}C^{\infty }(]0,T],\tilde{B}_{\infty }^{\sigma ,\infty }).$
La Proposition \ref{prop3} termine alors la preuve. $\blacksquare $

\subsection{Démonstration de la Proposition \protect\ref{prop1}}

Pour démontrer cette proposition, nous allons suivre une approche inspirée
du travail \cite{Che99} de J.-Y. Chemin. Nous décomposons la démonstration
en deux étapes:

\subsubsection{Première étape}

Soit $u_{0}\in L^{q}(\mathbb{R}^{d})$ et soit $u\in L_{T}^{\frac{2}{1-r}%
}\left( B_{\infty }^{-r,\infty }\right) \cap L_{T}^{p}L_{x}^{q}$ une
solution des équations de Navier-Stokes associées à la donnée initiale $%
u_{0}.$ On se propose de montrer qu'il existe $T_{0}\in ]0,T]$ tel que $u\in 
\mathbf{\tilde{L}}_{T_{0}}^{\frac{2}{1+r}}\left( B_{q}^{1+r,\infty }\right)
. $ Pour ce faire, les lemmes suivants nous seront fort utiles.

\begin{lemma}
\label{lemm1} Soient $\delta \in ]0,T],\rho \in \lbrack \frac{2}{1+r},\infty
\lbrack ,m\in \lbrack 1,+\infty ]$ et $\sigma >r.$ Alors l'opérateur liné%
aire $\mathbb{L}_{u}$ défini par:%
\begin{equation}
\mathbb{L}_{u}(f)=\sum_{k=1}^{2}\mathbb{L}_{oss}\left( \Pi _{k}(u,f)\right)
\label{lu}
\end{equation}%
est continu sur l'espace $\mathbf{\tilde{L}}_{\delta }^{\rho }\left(
B_{m}^{\sigma ,\infty }\right) $ et sa norme est majorée par $C\left\Vert
u\right\Vert _{L_{\delta }^{\frac{2}{1-r}}\left( B_{\infty }^{-r,\infty
}\right) }$ où $C$ est une constante indépendante de $\delta .$
\end{lemma}

\begin{lemma}
\label{lemm2} On pose $\omega =\mathbf{B}(u,u)$ et $\omega _{0}=\mathbb{L}%
_{u}(e^{t\Delta }u_{0})$ où $\mathbb{L}_{u}$ est l'opérateur défini par (\ref%
{lu}). Alors

\begin{enumerate}
\item $\omega \in \mathbf{\tilde{L}}_{T}^{\frac{p}{2}}\left(
B_{q/2}^{1,\infty }\right) $

\item $\omega _{0}\in \mathbf{\tilde{L}}_{T}^{\frac{2}{1+r}}\left(
B_{q}^{1+r,\infty }\right) $

\item $\omega _{0}\in \mathbf{\tilde{L}}_{T}^{\frac{p}{2}}\left( B_{q/2}^{1+%
\frac{2}{p},\infty }\right) .$
\end{enumerate}
\end{lemma}

\begin{lemma}
\label{lemm3} Soient $X_{1}$ et $X_{2}$ deux espaces de Banach et $f$ une
application définie sur $X_{1}+X_{2}$. Si les applications $%
f:X_{1}\rightarrow X_{1}$ et $f:X_{2}\rightarrow X_{2}$ sont contractantes
et si $z$ est un point fixe de $f$ dans $X_{1}$ alors $z\in X_{2}.$
\end{lemma}

Admettons un instant ces lemmes et montrons qu'il existe $T_{0}\in ]0,T]$
tel que $u$ soit dans l'espace $\mathbf{\tilde{L}}_{T_{0}}^{\frac{2}{1+r}%
}\left( B_{q}^{1+r,\infty }\right) .$ On pose $\omega =\mathbf{B}(u,u)$ et $%
\omega _{0}=\mathbb{L}_{u}(e^{t\Delta }u_{0})$. On considère la dé%
composition suivante de $\omega :$%
\begin{equation*}
\omega =\omega _{0}+\mathbb{L}_{u}(\omega )\equiv \digamma _{u}(\omega ).
\end{equation*}%
Les lemmes \ref{lemm1} et \ref{lemm2} assurent que, pour $T_{0}>0$ assez
petit (de sorte que $\left\Vert u\right\Vert _{L_{T_{0}}^{\frac{2}{1-r}%
}\left( B_{\infty }^{-r,\infty }\right) }$ soit inférieur à une constante
absolue $\varepsilon $ dépendant seulement de $r,p,q$), l'application affine 
$\digamma _{u}$\ est contractante sur les espaces de Banach $\mathbf{\tilde{L%
}}_{T_{0}}^{\frac{2}{1+r}}\left( B_{q}^{1+r,\infty }\right) $ et $\mathbf{%
\tilde{L}}_{T_{0}}^{\frac{p}{2}}\left( B_{q/2}^{1,\infty }\right) $ (c'est
ici qu'intervient l'hypothèse $p\geq \frac{4}{1+r}).$ Or, d'après le Lemme %
\ref{lemm2} et par construction de l'application $\digamma _{u},$ $\omega $
est un point fixe de cette application dans l'espace $\mathbf{\tilde{L}}%
_{T_{0}}^{\frac{2}{1+r}}\left( B_{q}^{1+r,\infty }\right) .$ Alors, le Lemme %
\ref{lemm3} implique $\omega \in \mathbf{\tilde{L}}_{T_{0}}^{\frac{2}{1+r}%
}\left( B_{q}^{1+r,\infty }\right) .$ Ce qui nous permet de conclure puisque 
$e^{t\Delta }u_{0}$ appartient aussi à $\mathbf{\tilde{L}}_{T_{0}}^{\frac{2}{%
1+r}}\left( B_{q}^{1+r,\infty }\right) $ (voir la démonstration du Lemme \ref%
{lemm2}). $\blacksquare $\medskip

Montrons maintenant les lemmes \ref{lemm1},\ref{lemm2} et \ref{lemm3}%
.\medskip

\noindent \textbf{Démonstration du lemme \ref{lemm1}:} Elle est une consé%
quence immédiate de l'injection continue de $L_{T}^{\frac{2}{1-r}}\left(
B_{\infty }^{-r,\infty }\right) $ dans $\mathbf{\tilde{L}}_{T}^{\frac{2}{1-r}%
}\left( B_{\infty }^{-r,\infty }\right) ,$ des lois du paraproduit de Bony
dans les espaces de Chemin-Lerner (Proposition \ref{pro3}) et de l'effet ré%
gularisant de l'opérateur intégral d'Oseen (Proposition \ref{pro5}). $%
\blacksquare \bigskip $

\noindent \textbf{Démonstration du lemme \ref{lemm2}:} Le premier point se dé%
montre aisément en se servant du fait que $u^{2}\in L_{T}^{p/2}L_{x}^{q/2}$
et en utilisant ensuite l'injection de cet espace dans $\mathbf{\tilde{L}}%
_{T}^{\frac{p}{2}}\left( B_{q/2}^{0,\infty }\right) $ et l'effet ré%
gularisant de l'opérateur intégral d'Oseen (Proposition \ref{pro5}). Pour
montrer le second point, il suffit d'utiliser l'injection de $L^{q}(\mathbb{R%
}^{d})$ dans $\tilde{B}_{q}^{0,\infty }$ qui implique, grâce à la
Proposition \ref{pro4}, que la \textit{tendence} $U_{0}\equiv e^{t\Delta
}u_{0}$ appartienne à $\mathbf{\tilde{L}}_{T}^{\frac{2}{1+r}}\left(
B_{q}^{1+r,\infty }\right) $ (et à $\mathbf{\tilde{L}}_{T}^{p}\left( B_{q}^{%
\frac{2}{p},\infty }\right) $ aussi), puis d'appliquer le lemme précédent.
Enfin, le dernier point résulte du fait que $e^{t\Delta }u_{0}\in \mathbf{%
\tilde{L}}_{T}^{p}\left( B_{q}^{\frac{2}{p},\infty }\right) ,$ de la
continuité des opérateurs du paraproduit de Bony 
\begin{equation*}
\Pi _{k}:L_{T}^{p}L_{x}^{q}\times \mathbf{\tilde{L}}_{T}^{p}\left( B_{q}^{%
\frac{2}{p},\infty }\right) \rightarrow \mathbf{\tilde{L}}_{T}^{\frac{P}{2}%
}\left( B_{\frac{q}{2}}^{\frac{2}{p},\infty }\right)
\end{equation*}%
et de l'effet régularisant de l'opérateur intégral d'Oseen. $\blacksquare
\bigskip $

\noindent \textbf{Démonstration du lemme \ref{lemm3}: }On munit l'espace $%
X\equiv X_{1}\cap X_{2}$ de la norme naturelle $\left\Vert x\right\Vert
=\left\Vert x\right\Vert _{X_{1}}+\left\Vert x\right\Vert _{X_{2}}.$ Il est
clair que $X$ est un espace de Banach et que $f$ est contractante sur cet
espace$,$ d'où l'existence d'un point fixe $z^{\prime }$ de $f$ dans $X.$
Ensuite l'unicité du point fixe de $f$ dans l'espace $X_{1}$ implique que $%
z=z^{\prime }$. Ce qui termine la preuve. $\blacksquare $

\subsubsection{Deuxième étape}

Soient $u_{1},u_{2}\in L_{T}^{\frac{2}{1-r}}\left( B_{\infty }^{-r,\infty
}\right) \cap L_{T}^{p}L_{x}^{q}$ deux solutions des équations de
Navier-Stokes associées à la même donnée initiale $u_{0}\in L^{q}(\mathbb{R}%
^{d}).$ D'après l'étape précédente, il existe $T_{0}\in ]0,T]$ tel que 
\begin{equation*}
u_{1},u_{2}\in \mathcal{Z}_{T_{0}}\equiv \mathbf{\tilde{L}}_{T_{0}}^{\frac{2%
}{1+r}}\left( B_{q}^{1+r,\infty }\right) \cap \mathbf{\tilde{L}}_{T_{0}}^{%
\frac{2}{1-r}}\left( B_{\infty }^{-r,\infty }\right) .
\end{equation*}%
Soit $\delta \in ]0,T_{0}]$ un réel arbitraire à choisir dans la suite$.$
Une application directe des lois du paraproduit de Bony dans les espaces de
Chemin-Lerner (Proposition \ref{pro3}) et de l'effet régularisant de l'opé%
rateur d'Oseen (Proposition \ref{pro5}), montre que l'application bilinéaire 
$\mathbf{B}$ est continue de $\mathcal{Z}_{\delta }\times \mathcal{Z}%
_{\delta }$ dans l'espace $\mathbf{\tilde{L}}_{\delta }^{\frac{2}{1+r}%
}\left( B_{q}^{1+r,\infty }\right) \cap \mathbf{\tilde{L}}_{\delta }^{\frac{2%
}{1-r}}\left( B_{q}^{1-r,\infty }\right) $ et que sa norme est majorée par
une constante $C$ indépendante de $\delta .$ Ensuite, comme $q\geq d$ alors
les injections de Bernstein (Proposition \ref{pro1}) assurent que $\mathbf{%
\tilde{L}}_{\delta }^{\frac{2}{1-r}}\left( B_{q}^{1-r,\infty }\right) $
s'injecte par continuité dans $\mathbf{\tilde{L}}_{\delta }^{\frac{2}{1-r}%
}\left( B_{\infty }^{-r,\infty }\right) $ et que la norme de cette injection
est indépendante de $\delta .$ En conclusion, l'application bilinéaire $%
\mathbf{B}:\mathcal{Z}_{\delta }\times \mathcal{Z}_{\delta }\rightarrow 
\mathcal{Z}_{\delta }$ est continue et sa norme est majorée par une
constante $C$ indépendante de $\delta .$ On a donc 
\begin{equation*}
\left\Vert u_{1}-u_{2}\right\Vert _{\mathcal{Z}_{\delta }}\leq C\left(
\left\Vert u_{1}\right\Vert _{\mathcal{Z}_{\delta }}+\left\Vert
u_{1}\right\Vert _{\mathcal{Z}_{\delta }}\right) \left\Vert
u_{1}-u_{2}\right\Vert _{\mathcal{Z}_{\delta }}.
\end{equation*}%
Ce qui implique, en prenant $\delta $ assez petit, que $u_{1}=u_{2}$ sur $%
[0,\delta ].$ $\blacksquare $

\subsection{Démonstration de la Proposition \protect\ref{prop2}}

Le premier point est classique (voir par exemple \cite{Can}, \cite{Lem02} et 
\cite{Mey}). Les démonstrations des deux autres points s'appuient
principalement sur le lemme élémentaire suivant où l'on utilise la notation
suivante:

\begin{notation}
Soient $T$ et $\mu $ deux réels positifs. On désigne par $L_{\mu ,T}^{\infty
}$ l'espace des fonctions $f:[0,T[\times \mathbb{R}^{d}\rightarrow \mathbb{R}%
^{d}$ telles que 
\begin{equation*}
\left\Vert f\right\Vert _{L_{\mu ,T}^{\infty }}\equiv \sup_{0<s<T}s^{\mu
/2}\left\Vert f(s)\right\Vert _{\infty }<\infty .
\end{equation*}
\end{notation}

\begin{lemma}
\label{lemme4} Soient $r\in ]0,1[$ et $T>0.$ Alors l'application bilinéaire $%
\mathbf{B}$ est continue de $L_{1,T}^{\infty }\times L_{r,T}^{\infty }$
(respectivement $L_{r,T}^{\infty }\times L_{r,T}^{\infty })$ dans $%
L_{r,T}^{\infty }$ et sa norme est majorée par $C_{r,d}$ (respectivement $%
C_{r,d}T^{\frac{1-r}{2}}$) où $C_{r,d}$ est une constante qui ne dépend que
de $r$ et de $d.$
\end{lemma}

La démonstration de ce lemme est immédiate. Il suffit de remarquer que l'opé%
rateur $e^{(t-s)\Delta }\mathbb{P}\nabla $ est un opérateur de convolution
avec une fonction intégrable dont la norme $L^{1}\left( \mathbb{R}%
^{d}\right) $ est de l'ordre de $(t-s)^{-1/2}$ (voir par exemple \cite{Can}, 
\cite{Lem02} et \cite{Mey}).\medskip

Revenons à la démonstration de deux dernières assertions de notre
proposition.\medskip

Il est bien connu (voir par exemple \cite{Lem02}) qu'il existe $T_{0}\in
]0,T]$ tel que la solution $u$ soit la limite dans l'espace $X_{T_{0}}\equiv
L_{T_{0}}^{\infty }L_{x}^{q}\cap L_{1,T_{0}}^{\infty }$ de la suite $%
(u_{(n)})_{n}$ définie par la relation de récurrence: 
\begin{equation*}
u_{(0)}=e^{t\Delta }u_{0};~\forall n\in \mathbb{N},~u_{(n+1)}=e^{t\Delta
}u_{0}+\mathbf{B}(u_{(n)},u_{(n)}),
\end{equation*}%
et que la série de terme général $\sigma _{n}\equiv \left\Vert
u_{(n+1)}-u_{(n)}\right\Vert _{X_{T_{0}}}$ soit convergente. Montrons que la
suite $(u_{(n)})_{n}$ est de Cauchy dans l'espace $L_{r,T_{0}}^{\infty }.$
Tout d'abord, comme $u_{0}\in B_{\infty }^{-r,\infty },$ alors la caracté%
risation des espaces de Besov à l'aide du noyau de la chaleur assure que $%
u_{(0)}\in L_{r,T_{0}}^{\infty }.$ Ensuite, en procédant par récurrence et
en utilisant le lemme précédent , on montre aisément que la suite $%
(u_{(n)})_{n}$ est dans $L_{r,T_{0}}^{\infty }$ et vérifie l'inégalité
suivante: 
\begin{equation*}
\left\Vert u_{(n+1)}-u_{(n)}\right\Vert _{L_{r,T_{0}}^{\infty }}\leq
C_{r,d}\sigma _{n}\left( \left\Vert u_{(n)}\right\Vert _{L_{r,T_{0}}^{\infty
}}+\left\Vert u_{(n-1)}\right\Vert _{L_{r,T_{0}}^{\infty }}\right)
\end{equation*}%
Ce qui implique (voir \cite{FLZZ} et \cite{Lem08}) que la suite $%
(u_{(n)})_{n}$ est de Cauchy dans $L_{r,T_{0}}^{\infty }.$ Rappelons à ce
stade que $L_{r,T_{0}}^{\infty }$ est un espace de Banach, qu'il s'injecte
par continuité dans $L_{1,T_{0}}^{\infty }$ et que la suite $(u_{(n)})_{n}$
converge vers $u$ dans $L_{1,T_{0}}^{\infty },$ il vient que $u\in
L_{r,T_{0}}^{\infty }$. Montrons maintenant que $u\in L_{T_{0}}^{\infty
}(B_{\infty }^{-r,\infty })$ ce qui achèvera la preuve de la deuxième
assertion grâce à (\ref{reg})$.$ Tout d'abord, on a bien d'après la
Proposition \ref{pro4}, $e^{t\Delta }u_{0}\in L_{T_{0}}^{\infty }(B_{\infty
}^{-r,\infty }).$ D'autre part, la Proposition \ref{pro2}, les inégalités de
Young et le fait que la norme $L^{1}(\mathbb{R}^{d})$ du noyau de l'opé%
rateur $e^{t\Delta }\mathbb{P}\nabla $ soit de l'ordre de $\frac{1}{\sqrt{t}}
$ impliquent que, pour tout $t\in ]0,T_{0}],$ 
\begin{eqnarray*}
\left\Vert \mathbf{B}(u,u)(t)\right\Vert _{B_{\infty }^{-r,\infty }}
&=&\sup_{0<\theta \leq 1}\theta ^{r/2}\left\Vert e^{\theta \Delta }\mathbf{B}%
(u,u)(t)\right\Vert _{\infty } \\
&\lesssim &\sup_{0<\theta \leq 1}\theta ^{r/2}\int_{0}^{t}\frac{1}{%
s^{(r+1)/2}\sqrt{t+\theta -s}}ds\left\Vert u\right\Vert
_{L_{r,T_{0}}^{\infty }}\left\Vert u\right\Vert _{L_{1,T_{0}}^{\infty }} \\
&\lesssim &\left\Vert u\right\Vert _{L_{r,T_{0}}^{\infty }}\left\Vert
u\right\Vert _{L_{1,T_{0}}^{\infty }}.
\end{eqnarray*}%
Ce qui conclut la preuve du second point.

Passons enfin à la démonstration de la dernière assertion de notre
proposition. Supposons que le temps maximal $T^{\ast }$ soit fini. Soit $%
r\in ]0,1[$ et soit $t_{0}$ un réel quelconque de l'intervalle $I_{\ast
}\equiv ]\max (0;T^{\ast }-1),T^{\ast }[$. D'après la Remarque \ref{mild},
la fonction $u_{t_{0}}(.)\equiv u(.+t_{0})$ est une solution mild sur $%
[0,\delta _{0}\equiv T^{\ast }-t_{0}[$ des équations de Navier-Stokes avec
donnée initiale $u(t_{0}).$ Alors, pour tout $t\in \lbrack 0,\delta _{0}[,$%
\begin{equation*}
u_{t_{0}}(t)=e^{t\Delta }u(t_{0})+\mathbf{B}(u_{t_{0}},u_{t_{0}})(t).
\end{equation*}%
Ainsi, en utilisant la Proposition \ref{pro2} et en appliquant le lemme précé%
dent , on en déduit qu'il existe une constante $C>0$ qui ne dépend que de $r$
et $d$ telle que, pour tout $t\in \lbrack 0,\delta _{0}[,$%
\begin{equation*}
t^{r/2}\left\Vert u_{t_{0}}(t)\right\Vert _{\infty }\leq C\left( \left\Vert
u(t_{0})\right\Vert _{B_{\infty }^{-r,\infty }}+t^{\frac{1-r}{2}}\left(
\sup_{0\leq s\leq t}s^{r/2}\left\Vert u_{t_{0}}(s)\right\Vert _{\infty
}\right) ^{2}\right) .
\end{equation*}%
Posons $f(t)\equiv \sup_{0\leq s\leq t}s^{r/2}\left\Vert
u_{t_{0}}(s)\right\Vert _{\infty }$; il vient,%
\begin{equation*}
\forall t\in \lbrack 0,\delta _{0}[,~f(t)\leq C\left( \left\Vert
u(t_{0})\right\Vert _{B_{\infty }^{-r,\infty }}+\left( T^{\ast
}-t_{0}\right) ^{\frac{1-r}{2}}f^{2}(t)\right) .
\end{equation*}%
Rappelons maintenant qu'il est bien connu (\cite{Gig}, \cite{Lem02}, \cite%
{May03}) que $\left\Vert u(t)\right\Vert _{\infty }\rightarrow +\infty $
lorsque $t\rightarrow T^{\ast },$ ce qui implique que $f(t)\rightarrow
+\infty $ lorsque $t\rightarrow \delta _{0},$ et utilisons ensuite le lemme %
\ref{lemme5} ci-dessous, il vient 
\begin{equation*}
\left\Vert u(t_{0})\right\Vert _{B_{\infty }^{-r,\infty }}\left( T^{\ast
}-t_{0}\right) ^{\frac{1-r}{2}}\geq \varepsilon _{r,d}\equiv \frac{1}{4C^{2}}%
.
\end{equation*}%
Ce qui termine la preuve de la Proposition \ref{prop2}.

\begin{lemma}
\label{lemme5} soient $a<b$ deux réels et $f:[a,b[\rightarrow \mathbb{R}$
une application continue. On suppose qu'il existe deux réels $A$ et $B>0$
tels que: $4AB<1,f(0)\leq 2A$ et $\forall t\in \lbrack a,b[,~f(t)\leq
A+Bf^{2}(s).$ Alors $\forall t\in \lbrack a,b[,~f(t)\leq 2A.$
\end{lemma}

La démonstration de ce lemme est simple. Il suffit de remarquer qu'en vertu
de l'hypothèse $4AB<1,$ $f$ ne peut pas prendre la valeur $2A$ et
d'appliquer ensuite le théorème des valeurs intermédiaires. $\blacksquare $

\subsection{Démonstration de la Proposition \protect\ref{prop3}}

La démonstration de cette proposition s'inspire de l'article \cite{Lem07} de
P.-G. Lemarié-Rieusset. Soit $(t_{n})_{n}\in ]0,T/2]$ une suite qui tend
vers $0.$ On considère la suite de fonctions $\left( u_{n}\right) _{n}$ dé%
finie sur $[0,T/2[$ par $u_{n}(t)=u(t+t_{n}).$ Il s'agit de montrer que $%
\sup_{0<t<\delta }\sqrt{t}\left\Vert u_{n}(t)\right\Vert _{\infty }$ tend
vers $0$ uniformément par rapport à $n$ lorsque $\delta $ tend vers $0.$
Tout d'abord, pour alléger l'écriture, nous introduisons les notations
suivantes: \ 
\begin{eqnarray*}
h_{n}(\mu ,\delta ) &\equiv &\sup_{0<t\leq \delta }t^{\frac{\mu +1}{2}%
}\left\Vert u_{n}(t)\right\Vert _{B_{\infty }^{\mu ,\infty }}, \\
\Theta (\delta ) &=&\sup_{0<t_{0}<\frac{T}{2}}\left\Vert u\right\Vert _{L^{%
\frac{2}{1-r}}([t_{0},t_{0}+\delta ],B_{\infty }^{-r,\infty })}.
\end{eqnarray*}%
Soient $\sigma \in ]r,1[$ un réel fixe et $\delta _{0}\in ]0,T/2[$ à choisir
ultérieurement. Soient $n\in \mathbb{N},\delta \in ]0,\delta _{0}]$ et $t\in
]0,\delta ].$ On désigne par $a=a(n,t)$ un réel appartenant à $[\frac{t}{4},%
\frac{t}{2}]$ tel que 
\begin{equation*}
\left\Vert u_{n}(a)\right\Vert _{B_{\infty }^{-r,\infty }}=\inf_{s\in
\lbrack \frac{t}{4},\frac{t}{2}]}\left\Vert u_{n}(s)\right\Vert _{B_{\infty
}^{-r,\infty }}.
\end{equation*}%
Comme $u_{n}$ est une solution mild des équations de Navier-Stokes, alors
d'après la Remarque \ref{mild}, 
\begin{eqnarray}
u_{n}(t) &=&e^{(t-a)\Delta }u_{n}(a)-\int_{a}^{t}e^{(t-s)\Delta }\mathbb{P}%
\nabla \left( u_{n}\otimes u_{n}\right) ds  \label{ee1} \\
&\equiv &I_{n}(t)+J_{n}(t).  \label{ee2}
\end{eqnarray}%
Nous allons estimer les normes de $I_{n}(t)$ et $J_{n}(t)$ dans l'espace de
Besov $B_{\infty }^{\sigma ,\infty }.$ D'après la première assertion de la
Proposition \ref{pro4} et la définition de $a=a(n,t),$ on a 
\begin{eqnarray}
\left\Vert I_{n}(t)\right\Vert _{B_{\infty }^{\sigma ,\infty }} &\lesssim
&(t-a)^{-\frac{r+\sigma }{2}}\left\Vert u_{n}(a)\right\Vert _{B_{\infty
}^{-r,\infty }}  \notag \\
&\lesssim &t^{-\frac{1+\sigma }{2}}\left\Vert u_{n}\right\Vert _{L^{\frac{2}{%
1-r}}([t/4,t/2],B_{\infty }^{-r,\infty })}  \notag \\
&\lesssim &t^{-\frac{1+\sigma }{2}}\Theta (\delta ).  \label{e1}
\end{eqnarray}%
D'autre part, en remaquant que $J_{n}(t)=\mathbb{L}_{oss}(1_{[a,t]}u_{n}%
\otimes 1_{[a,t]}u_{n})(t)$ et en utilisant la continuité des opérateurs du
paraproduit de Bony $\Pi _{k}$ de $\tilde{L}_{T}^{\frac{2}{1-r}}(B_{\infty
}^{-r,\infty })\times L_{T}^{\infty }\left( B_{\infty }^{\sigma ,\infty
}\right) $ dans $\tilde{L}_{T}^{\frac{2}{1-r}}(B_{\infty }^{\sigma -r,\infty
})$ (Proposition \ref{pro3}) et la continuité de l'opérateur $\mathbb{L}%
_{Oss}$ de $\tilde{L}_{T}^{\frac{2}{1-r}}(B_{\infty }^{\sigma -r,\infty })$
dans $L_{T}^{\infty }\left( B_{\infty }^{\sigma ,\infty }\right) $
(Proposition \ref{pro5}), on déduit aisément que 
\begin{eqnarray}
\left\Vert J_{n}(t)\right\Vert _{B_{\infty }^{\sigma ,\infty }} &\lesssim
&\left\Vert u_{n}\right\Vert _{\tilde{L}^{\frac{2}{1-r}}([a,t],B_{\infty
}^{-r,\infty })}\left\Vert u_{n}\right\Vert _{L^{\infty }([a,t],B_{\infty
}^{\sigma ,\infty })}  \notag \\
&\lesssim &\left\Vert u_{n}\right\Vert _{L^{\frac{2}{1-r}}([a,t],B_{\infty
}^{-r,\infty })}\sup_{a\leq s\leq t}\left\Vert u_{n}(s)\right\Vert
_{B_{\infty }^{\sigma ,\infty }}  \notag \\
&\lesssim &t^{-\frac{1+\sigma }{2}}\Theta (\delta )h_{n}(\sigma ,\delta ) 
\notag \\
&\lesssim &t^{-\frac{1+\sigma }{2}}\Theta (\delta _{0})h_{n}(\sigma ,\delta
).  \label{e2}
\end{eqnarray}%
Les inégalités (\ref{e1}) et (\ref{e2}) impliquent qu'il existe une
constante $C_{1}>0$ indépendente de $t,\delta $ et $n$ telle que%
\begin{equation}
h_{n}(\sigma ,\delta )\leq C_{1}\Theta (\delta )+C_{1}\Theta (\delta
_{0})h_{n}(\sigma ,\delta ).  \label{e3}
\end{equation}%
En choisissant alors $\delta _{0}$ assez petit de sorte que $\Theta (\delta
_{0})$ soit inférieure à $\frac{1}{2C_{1}}$ (ce qui est posssible puisque $%
\Theta (\delta _{0})\rightarrow 0$ lorsque $\delta _{0}\rightarrow 0$), l'iné%
galité précédente nous donne 
\begin{equation}
h_{n}(\sigma ,\delta )\leq 2C_{1}\Theta (\delta ).  \label{e4}
\end{equation}%
Nous retournons maintenant aux égalités (\ref{ee1}) et (\ref{ee2}) pour
estimer cette fois les normes de $I_{n}(t)$ et $J_{n}(t)$ dans l'espace de
Besov $B_{\infty }^{-r,\infty }.$ A nouveau, d'après la première assertion
de la Proposition \ref{pro4} et la définition de $a=a(n,t),$ on a 
\begin{eqnarray}
\left\Vert I_{n}(t)\right\Vert _{B_{\infty }^{-r,\infty }} &\lesssim
&\left\Vert u_{n}(a)\right\Vert _{B_{\infty }^{-r,\infty }}  \notag \\
&\lesssim &t^{-\frac{1-r}{2}}\Theta (\delta ).  \label{e5}
\end{eqnarray}%
D'autre part, en utilisant la continuité des opérateurs de Bony $\Pi _{k}$
de $B_{\infty }^{-r,\infty }\times B_{\infty }^{\sigma ,\infty }$ dans $%
B_{\infty }^{\sigma -r,\infty }$, l'action de l'opérateur pseudo-diffé%
rentiel $\mathbb{P}\nabla $ sur les espaces de Besov (Proposition \ref{pro1}%
) et la Proposition \ref{pro4} on en déduit aisément les estimations
suivantes:%
\begin{eqnarray}
\left\Vert J_{n}(t)\right\Vert _{B_{\infty }^{-r,\infty }} &\lesssim
&\int_{a}^{t}\frac{1}{(t-s)^{\frac{1-\sigma }{2}}}\left\Vert \mathbb{P}%
\nabla \left( u_{n}\times u_{n}\right) \right\Vert _{B_{\infty }^{\sigma
-r-1,\infty }}ds  \notag \\
&\lesssim &t^{\frac{1+\sigma }{2}}\sup_{\frac{t}{4}<s<t}\left\Vert
u_{n}(s)\right\Vert _{B_{\infty }^{-r,\infty }}\sup_{\frac{t}{4}%
<s<t}\left\Vert u_{n}(s)\right\Vert _{B_{\infty }^{\sigma ,\infty }}  \notag
\\
&\lesssim &t^{-\frac{1-r}{2}}h_{n}(-r,\delta )h_{n}(\sigma ,\delta )  \notag
\\
&\lesssim &t^{-\frac{1-r}{2}}h_{n}(-r,\delta )\Theta (\delta _{0})
\label{e6}
\end{eqnarray}%
où nous avons utilisé (\ref{e4}) dans le dernier passage.

Alors, en combinant les inégalités (\ref{ee1}) et (\ref{ee2}), on en déduit
qu'il existe une constante $C_{2}>0$ indépendente de $t,\delta $ et $n$
telle que%
\begin{equation*}
h_{n}(-r,\delta )\leq C_{2}\Theta (\delta )+C_{2}\Theta (\delta
_{0})h_{n}(-r,\delta )
\end{equation*}%
Ainsi, pour $\delta _{0}$ assez petit, on a 
\begin{equation}
h_{n}(-r,\delta )\leq 2C_{2}\Theta (\delta )  \label{e40}
\end{equation}%
En résumé, en se servant, comme dans \cite{Lem07} et \cite{LM}, de l'inégalit%
é classique d'interpolation 
\begin{equation*}
\left\Vert f\right\Vert _{\infty }\leq \left( \left\Vert f\right\Vert
_{B_{\infty }^{-r,\infty }}\right) ^{\frac{\sigma }{r+s}}\left( \left\Vert
f\right\Vert _{B_{\infty }^{\sigma ,\infty }}\right) ^{\frac{r}{r+\sigma }},
\end{equation*}%
les inégalités (\ref{e4}) et (\ref{e40}) impliquent qu'il existe deux
constantes $C>0$ et $\delta _{0}\in ]0,T/2]$ indépendantes de $n$ telles que
pour tout $\delta \in ]0,\delta _{0}]$ on ait 
\begin{equation*}
\sup_{0<t<\delta }\sqrt{t}\left\Vert u_{n}(t)\right\Vert _{\infty }\leq
C\Theta (\delta ).
\end{equation*}%
Ce qui termine la preuve de la Proposition \ref{prop3}$.\blacksquare $

\bigskip 
\providecommand{\bysame}{\leavevmode\hbox
to3em{\hrulefill}\thinspace}



\providecommand{\bysame}{\leavevmode\hbox
to3em{\hrulefill}\thinspace} \providecommand{\MR}{\relax\ifhmode\unskip%
\space\fi MR } 
\providecommand{\MRhref}[2]{  \href{http://www.ams.org/mathscinet-getitem?mr=#1}{#2}
} \providecommand{\href}[2]{#2}

\end{document}